\documentclass[10pt]{article}
\usepackage{fullpage}
\usepackage{cite} % for e.g. [1-3]
\usepackage{setspace}
\usepackage{graphicx}
\usepackage{delarray}
\usepackage{bbm}
\usepackage{amssymb}
\usepackage{amsmath}% needed for operatorwithlimits
\usepackage{amsthm}
\usepackage{algorithm}
\usepackage{algorithmic}
    \usepackage{sectsty}
    \sectionfont{\normalsize}

\theoremstyle{dotless}

\newtheorem{thm}{Theorem}

\newcommand{\be}{\begin{equation}}
\newcommand{\ee}{\end{equation}}
\newcommand{\bea}{\begin{eqnarray}}
\newcommand{\eea}{\end{eqnarray}}
\newcommand{\bfl}{\begin{flalign}}
\newcommand{\efl}{\end{flalign}}
\newcommand{\bfc}{\begin{figure}\begin{center}}
\newcommand{\efc}{\end{center}\end{figure}}

\newcommand{\dd}{\mathrm{d}}
\newcommand{\nin}{\noindent}

\title{Log-Quadratic Bounds for the Gaussian Q-function}

\author{Andrew Mastin, Patrick Jaillet\footnote{Both authors are with the Laboratory for Information and Decision Systems, Department of Electrical Engineering and Computer Science, Massachusetts Institute of Technology, Cambridge, MA 02139, USA {\tt \{mastin,jaillet\}@mit.edu}.
Supported by NSF grant 1029603. The first author is supported in part by a NSF graduate research fellowship.}}

\date{}
\begin{document}

\maketitle

\begin{abstract}
We present bounds of quadratic form for the logarithm of the Gaussian Q-function. We also show an analytical method for deriving log-quadratic approximations of the Q-function and give an approximation with absolute error less than $10^{-3}$.

\end{abstract}

\section{Introduction}
Approximations and bounds for the Gaussian Q-function have been studied extensively, motivated by the fact that the Q-function is not an elementary function \cite{mills26, gordon41, birnbaum42, tate53, boyd59, borj79, chiani03, deabreu09, chang11}\footnote{The Q-function is the tail probability for a standard normal random variable: \bea
Q(x) \triangleq \frac{1}{\sqrt{2\pi}} \int_x^\infty e^{-\frac{u^2}{2}} \dd u.
\eea}. There has been specific effort in deriving bounds that remain tractable under multiplication and exponentiation, especially those with a single exponential term. The classical example is the Chernoff bound \cite{chernoff52},
\be
Q(x) \le \frac{1}{2}e^{-\frac{x^2}{2}}, \quad x \ge 0.
\ee
In recent work, C{\^o}t{\'e} et al. \cite{cote12} showed strong Chernoff-type lower bounds for the Q-function of the form 
\be
Q(x) \ge \alpha e^{-\beta x^2}, \quad x \in \mathbb{R}.
\ee
Chernoff-type bounds are sufficient for many applications but may be loose for small values of $x$.  L{\'o}pez-Ben{\'\i}tez and Casadevall \cite{lopez11} considered Q-function approximations of the log-quadratic form
\bea
Q(x) &\approx& e^{-ax^2-bx-c}, \quad x \ge 0,
\eea
and determined the parameters $a,b,c$ using numerical fitting techniques for different ranges of $x$. We consider bounds and approximations of this form from an analytical perspective. The bounds are stated by the following theorem.

\begin{thm}
\bea
Q(x) &\ge& \frac{1}{2} e^{-\frac{1}{2}x^2 - \sqrt{\frac{2}{\pi}}x}, \quad x \ge 0, \\
Q(x) &\le& \frac{1}{2} e^{-\frac{1}{\pi}x^2 - \sqrt{\frac{2}{\pi}}x}, \quad x \ge 0.
\eea
\end{thm}
\nin We derive the following approximation, which has absolute error less than $10^{-3}$.
\bea
Q(x) &\approx& \frac{1}{2} e^{-0.374 x^2 - 0.777 x}, \quad x \ge 0.
\eea

\section{Bounds and Approximations}
We show the proof of Theorem 1 and then the derivation of the approximation. ~\\ ~\\
\textit{Proof of Theorem 1.} 
Let the standard normal density be given by
\be
\phi(x) \triangleq \frac{e^{-\frac{x^2}{2}}}{\sqrt{2\pi}}. 
\ee
Define
\bea
L(x) &\triangleq& \log(Q(x)), \\
\ell(x) &\triangleq& \frac{\partial L(x)}{\partial x} = -\frac{\phi(x)}{Q(x)}, \\
\ell'(x) &\triangleq& \frac{\partial^2 L(x)}{\partial^2 x} =\frac{\phi(x)Q(x)x-\phi(x)^2}{Q(x)^2} , \\
L_l(x) &\triangleq& -\frac{1}{2}x^2-\sqrt{\frac{2}{\pi}}x-\log(2), \\
\ell_l(x) &\triangleq& \frac{\partial L_l(x)}{\partial x} = -x-\sqrt{\frac{2}{\pi}}, \\
\ell_l'(x) &\triangleq& \frac{\partial^2 L_l(x)}{\partial x^2} = -1.
\eea
For the lower bound, we show that $L_l(x) \le L(x)$ for $x\ge0$. It is clear that $L_l(0) = L(0)$, so it is sufficient to show that $\ell_l(x) \le \ell(x)$ for $x \ge 0$. Again, we have $\ell_l(0) = \ell(0) = -\sqrt{\frac{2}{\pi}}$, so we are left to show that $\ell_l'(x) \le \ell'(x)$, or equivalently $\ell'(x) \ge -1$, for $x \ge 0$. The desired property 
\bea
\ell'(x) &=& \frac{\phi(x)Q(x)x-\phi(x)^2}{Q(x)^2} \ge -1
\eea
is equivalent to the statement
\bea
Q(x)^2+Q(x)\phi(x)x-\phi(x)^2 \ge 0.
\eea
The above inequality follows from a simple bound on Mill's ratio (e.g. \cite{sampford53}).

The upper bound is found with similar reasoning. It is sufficient to show that $\ell'(x) \le -\frac{2}{\pi}$, which is equivalent to the property
\bea
\frac{2}{\pi}Q(x)^2+Q(x)\phi(x)x-\phi(x)^2 \le 0.
\eea
An upper bound referred to in \cite{borj79} and \cite{dasgupta10} (as the Mitrinovi{\'c} inequality) is
\bea
Q(x) \le \frac{2\phi(x)}{x+\sqrt{x^2+\frac{8}{\pi}}}.
\eea
The substitution yields
\bea
\frac{2}{\pi}Q(x)^2+Q(x)\phi(x)x-\phi(x)^2 &\le& \frac{\left(\frac{8}{\pi}\right)\phi(x)^2}{\left(x+\sqrt{x^2+\frac{8}{\pi}}\right)^2} + \frac{2x\phi(x)^2}{x+\sqrt{x^2+\frac{8}{\pi}}} - \phi(x)^2 = 0.
\eea\qed
 ~\\
 \clearpage
\nin\textit{Q-function approximation}.
The approximation is of the form
\bea
Q(x) &\approx& e^{-ax^2-bx-c} ~\triangleq~ \hat{Q}(x).
\eea
for $x \ge 0$.
We use the function
\bea
\hat{L}(x) &\triangleq&  \log(\hat{Q}(x)) ~=~ -ax^2-bx-c
\eea
and select parameters to approximately minimize the mean squared error between derivatives of $L(x)$ and $\hat{L}(x)$ for $0 \le x \le t$. Define
\bea
\hat{\ell}(x,a,b) &\triangleq& \frac{\partial \hat{L}(x)}{\partial x} = -2ax-b, \\
\hat{\ell}'(x,a) &\triangleq& \frac{\partial^2 \hat{L}(x)}{\partial^2 x} = -2a.
\eea
The mean squared value for the difference of second derivatives of the two terms is proportional to
\bea
\int_0^t (\ell'(x)-\hat{\ell}'(x,a))^2 \dd x.
% &=& \int_0^t \left (2a -\frac{x \phi(x)}{Q(x)}+ \frac{\phi(x)^2}{Q(x)^2} \right )^2\dd x.
\eea
Differentiating with respect to $a$ and equating to zero gives an optimal value of $a$ for a given $t$:
\bea
\int_0^t \ell'(x) \dd x &=& \int_0^t \hat{\ell'}(x,a) \dd x, \\
\int_0^t\ell'(x) \dd x &=& \frac{\phi(0)}{Q(0)} - \frac{\phi(t)}{Q(t)} ~=~ \sqrt{\frac{2}{\pi}} - \frac{\phi(t)}{Q(t)}, \\
\int_0^t \hat{\ell'}(x,a) \dd x &=& -2at, \\
a^*(t) &=& -\frac{1}{2t}\left (\sqrt{\frac{2}{\pi}} - \frac{\phi(t)}{Q(t)} \right ).
\eea
The mean squared value for the difference of first derivatives is proportional to
\bea
\int_0^t (\ell(x)-\hat{\ell}(x,a,b))^2 \dd x.
\eea
Differentiating with respect to $b$, equating to zero, and substituting $a = a^*(t)$ gives the optimal value of $b$ as a function of $t$,
\bea
\int_0^t \ell(x) \dd x &=& \int_0^t \hat{\ell}(x,a,b) \dd x, \\
\int_0^t \ell(x)\dd x &=& \log(Q(t)) - \log(Q(0)) ~=~ \log \left(2Q(t) \right), \\
\int_0^t \hat{\ell}(x,a,b)\dd x &=& -at^2-bt, \\
b^*(t) &=& -a^*(t) t - \frac{1}{t} \log(2Q(t)), \\
b^*(t)&=& \frac{1}{2}\left (\sqrt{\frac{2}{\pi}} - \frac{\phi(t)}{Q(t)} \right ) - \frac{1}{t} \log(2Q(t)).
\eea
We simply use $c = \log(2)$ to match $Q(0) = \frac{1}{2}$. This gives the general approximation for a selected $t$
\bea
Q(x) \approx \frac{1}{2} \exp \left ({\frac{1}{2t}\left (\sqrt{\frac{2}{\pi}} - \frac{\phi(t)}{Q(t)} \right )x^2 + \left ( \frac{1}{t} \log(2Q(t)) - \frac{1}{2}\left (\sqrt{\frac{2}{\pi}} - \frac{\phi(t)}{Q(t)} \right ) \right )x} \right),\quad x \ge 0.
\eea
For a given $t$, an initial evaluation $Q(t)$ is required.

By making the parameter $t$ large, the approximation becomes increasingly more accurate for large values of $x$ and less accurate for small values of $x$. Selecting $t = 1.295$ approximately minimizes the maximum absolute error $|Q(x) - \hat{Q}(x)|$ for all $x \ge 0$; the corresponding error is $9.485\times 10^{-4}$. This gives the approximation
\be
Q(x) \approx \frac{1}{2} e^{-0.374 x^2 - 0.777 x}, \quad x \ge 0.
\ee

\qed

\bibliographystyle{jota}
\bibliography{q_bound}

\begin{thebibliography}{10}

\bibitem{mills26}
Mills, J.P.: Table of the ratio: area to bounding ordinate, for any portion of
  normal curve.
\newblock Biometrika pp. 395--400 (1926)

\bibitem{gordon41}
Gordon, R.D.: Values of {M}ills' ratio of area to bounding ordinate and of the
  normal probability integral for large values of the argument.
\newblock The Annals of Mathematical Statistics pp. 364--366 (1941)

\bibitem{birnbaum42}
Birnbaum, Z.: An inequality for {M}ill's ratio.
\newblock The Annals of Mathematical Statistics 13, 245--246 (1942)

\bibitem{tate53}
Tate, R.F.: On a double inequality of the normal distribution.
\newblock The Annals of Mathematical Statistics pp. 132--134 (1953)

\bibitem{boyd59}
Boyd, A.: Inequalities for {M}ills' ratio.
\newblock Rep. Statist. Appl. Res. Un. Japan. Sci. Engrs 6, 44--46 (1959)

\bibitem{borj79}
Borjesson, P., Sundberg, C.E.: Simple approximations of the error function
  {Q}(x) for communications applications.
\newblock Communications, IEEE Transactions on 27, 639--643 (1979)

\bibitem{chiani03}
Chiani, M., Dardari, D., Simon, M.K.: New exponential bounds and approximations
  for the computation of error probability in fading channels.
\newblock Wireless Communications, IEEE Transactions on 2, 840--845 (2003)

\bibitem{deabreu09}
de~Abreu, G.T.F.: Jensen-{C}otes upper and lower bounds on the {G}aussian
  {Q}-function and related functions.
\newblock Communications, IEEE Transactions on 57, 3328--3338 (2009)

\bibitem{chang11}
Chang, S.H., Cosman, P.C., Milstein, L.B.: Chernoff-type bounds for the
  {G}aussian error function.
\newblock Communications, IEEE Transactions on 59, 2939--2944 (2011)

\bibitem{chernoff52}
Chernoff, H.: A measure of asymptotic efficiency for tests of a hypothesis
  based on the sum of observations.
\newblock The Annals of Mathematical Statistics 23, 493--507 (1952)

\bibitem{cote12}
C{\^o}t{\'e}, F.D., Psaromiligkos, I.N., Gross, W.J.: A {C}hernoff-type lower
  bound for the {G}aussian {Q}-function.
\newblock arXiv preprint arXiv:1202.6483  (2012)

\bibitem{lopez11}
L{\'o}pez-Ben{\'\i}tez, M., Casadevall, F.: Versatile, accurate, and
  analytically tractable approximation for the {G}aussian {Q}-function.
\newblock Communications, IEEE Transactions on 59, 917--922 (2011)

\bibitem{sampford53}
Sampford, M.R.: Some inequalities on {M}ill's ratio and related functions.
\newblock The Annals of Mathematical Statistics pp. 130--132 (1953)

\bibitem{dasgupta10}
DasGupta, A.: Fundamentals of probability: a first course.
\newblock Springer Berlin (2010)

\end{thebibliography}

\end{document}